\documentclass[pdflatex,sn-mathphys-num]{sn-jnl}%

\usepackage{graphicx}%
\usepackage{multirow}%
\usepackage{amsmath,amssymb,amsfonts}%
\usepackage{amsthm}%
\usepackage{mathrsfs}%
\usepackage[title]{appendix}%
\usepackage{xcolor}%
\usepackage{textcomp}%
\usepackage{manyfoot}%
\usepackage{booktabs}%
\usepackage{algorithm}%
\usepackage{algorithmicx}%
\usepackage{algpseudocode}%
\usepackage{listings}%
\usepackage{lineno}
\usepackage{subcaption}
\usepackage{float}

\theoremstyle{thmstyleone}%
\theoremstyle{thmstyletwo}%
\theoremstyle{thmstylethree}%
\begin{document}

\title[]{Scalable Solution of the Stochastic Multi-path Traveling Salesman Problem via Neural Networks}


\author*[1]{\fnm{Xiaochen} \sur{Chou}}\email{xiaochen.chou@unimib.it}

\author[1]{\fnm{Ludovica} \sur{Di Marco}}\email{l.dimarco8@campus.unimib.it}

\author[1]{\fnm{Enza} \sur{Messina}}\email{enza.messina@unimib.it}

\affil*[1]{\orgdiv{Department of Informatics, Systems and Communication}, \orgname{University of Milano-Bicocca}, \orgaddress{\street{Viale Sarca 336}, \city{Milan}, \postcode{20125}, \country{Italy}}}
            

\abstract{The multi-path Traveling Salesman Problem with stochastic travel costs arises in hybrid vehicle routing applications designed for Smart City and City Logistics, where multiple paths exist between each pair of locations. Travel times along these paths are typically affected by real-time traffic conditions and therefore modeled as stochastic. The objective of the problem is to determine a Hamiltonian tour that minimizes the expected total travel cost under uncertainty.

In this work, we adopt a two-stage stochastic programming formulation. In the first stage, a predefined route specifying the sequence of locations to be visited is determined, while taking into consideration a second-stage recourse problem that selects the optimal path from the feasible set of alternative paths for each pair of locations, once real-time traffic conditions are realized. 
To reduce the computational burden imposed by the large number of scenarios required to capture travel time uncertainty, the innovation of this work is the integration of neural network-based surrogate models to approximate the expected value of the second-stage recourse problem. Different architectures and training strategies for the neural networks are proposed and analyzed, with performance evaluated in terms of computation time, solution quality, and generalization capability. Preliminary findings demonstrate the enhanced scalability and practical applicability of the approach for complex vehicle routing problems under uncertainty.}



\maketitle

\section{Introduction}\label{sec:intro}
The multi-path Traveling Salesman Problem (mpTSP) \cite{Tadei17, Perboli17} is a variant of the classical Traveling Salesman Problem (TSP), which aims to determine a tour with minimum cost that visits each location exactly once. In realistic urban mobility and logistics settings, travel costs between locations are often uncertain, therefore modeled with stochastic travel times \cite{Toriello14} or robust time intervals \cite{Montemanni07}. In the mpTSP, these costs are not represented on a single arc but by multiple alternative paths (e.g., different road corridors) connecting each pair of locations. Consequently, the problem is defined on a multipath network, where each location pair is linked by several candidate paths with distinct, typically traffic-dependent costs. The objective is then to identify a tour that minimizes the expected total cost over this stochastic, multipath network.

\begin{figure}[t]
    \centering
    \begin{subfigure}[b]{0.49\textwidth}
        \centering
        \includegraphics[width=\textwidth]{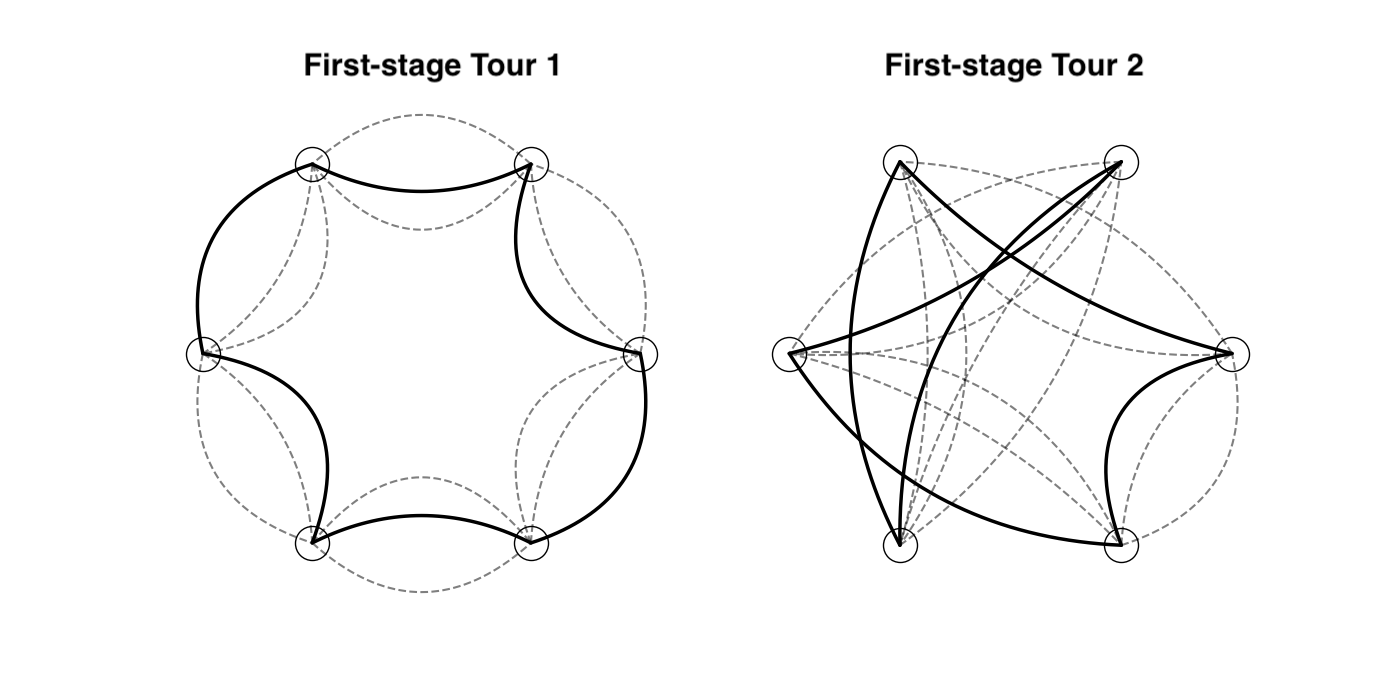}
        \caption{Two distinct first-stage tours with the corresponding second-stage path selections under a fixed scenario realization.}
        \label{fig:1a}
    \end{subfigure}
    \hfill
    \begin{subfigure}{0.49\textwidth}
        \centering
        \includegraphics[width=\textwidth]{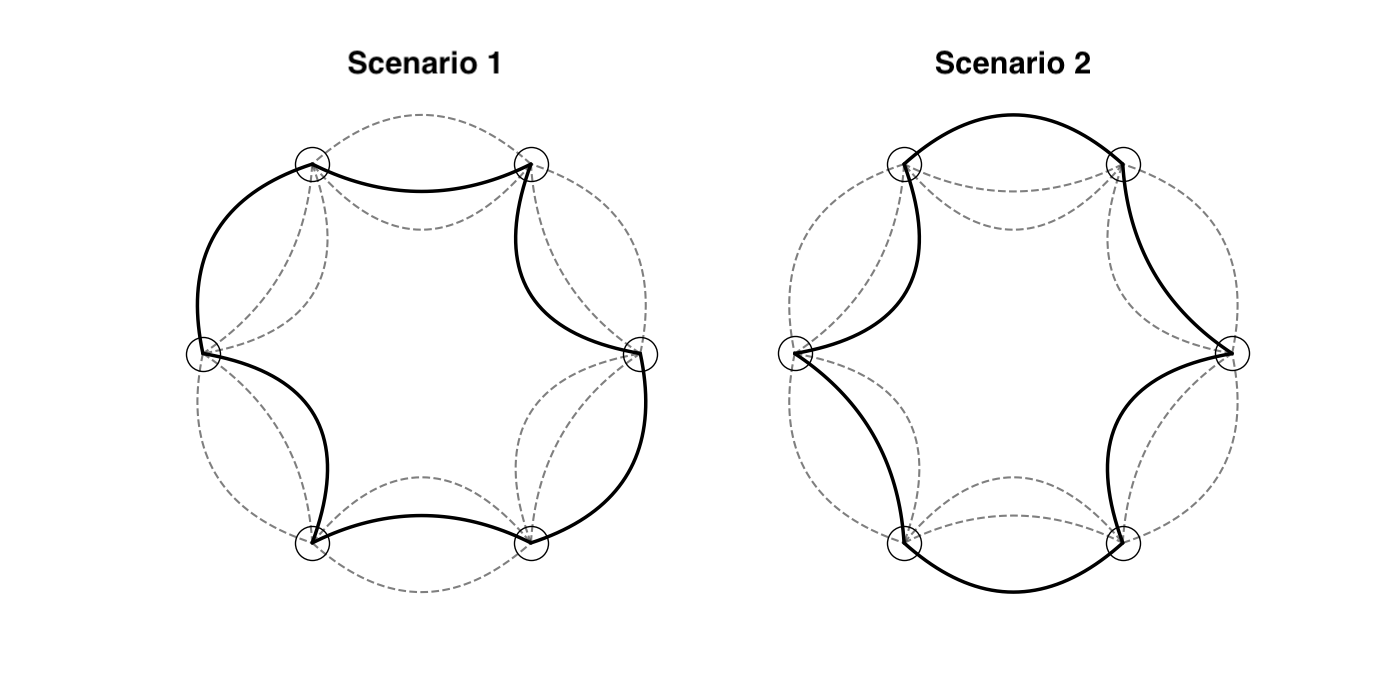}
        \caption{Second-stage path selections for a fixed first-stage tour across two different stochastic scenario realizations.}
        \label{fig:1b}
    \end{subfigure}
    \caption{Illustration of the mpTSP: first-stage and second-stage decisions}
    \label{fig:mpTSP}
\end{figure}

Stochastic programming (SP) \cite{SPbook} is a widely used framework for adapting models to stochastic settings. The mpTSP is formulated as a two-stage stochastic programming problem in \cite{Perboli17}: In the first stage, binary decision variables are used to select arcs, thereby determining a tour that visits the required locations in a specific order. The second-stage recourse problem selects the path to be traveled between each pair of locations in the tour, based on stochastic travel time realizations. We define a \emph{scenario} as a realization of traffic conditions during a particular time period, in which the travel times on all available paths between connected locations are known and deterministic. The overall objective is to identify the optimal first-stage tour while accounting for adaptive second-stage routing decisions across possible scenarios. Figure \ref{fig:mpTSP} illustrates the interaction between the two-stage decisions in the mpTSP. Figure \ref{fig:1a} presents the two different first-stage tours under a specific given scenario, along with their corresponding second-stage path selections. Figure \ref{fig:1b} illustrates the variability in second-stage path selections arising from two distinct scenario realizations under a fixed first-stage tour.

In routing problems where traffic conditions are considered, large scenario sets can be generated from sensor data to empirically approximate the underlying distribution of stochastic variability. However, incorporating all such data directly into the model makes it computationally intractable, as the complexity of the stochastic formulation grows rapidly with the increasing number of scenarios \cite{Perboli17, SPbook}. In mpTSP, this is primarily because the first-stage problem, even under deterministic settings, corresponds to a TSP, which is NP-hard \cite{Garey79}. In contrast, the second-stage (recourse) problem of the mpTSP is relatively easy to solve. The recourse problem of the mpTSP evaluates the performance of a given first-stage tour under specific scenario realizations. For a fixed first-stage solution and a known scenario set, this evaluation is computationally straightforward, since the shortest paths between locations can be readily identified for each scenario. Therefore, the computational challenge lies in determining the optimal first-stage decision while accounting for a large number of stochastic realizations.

To address this issue, surrogate models can be employed to approximate the recourse function value, providing efficient estimation of the expected performance of first-stage solutions. Two main modeling approaches have been proposed in the literature \cite{Dumouchelle22, Chou25}. The first approach integrates scenario information directly into the learning process by training neural networks that map first-stage decisions and scenario embeddings to recourse function values \cite{Dumouchelle22}. This design allows the surrogate to approximate the conditional response of the system under varying stochastic environments. The second approach, in contrast, focuses on learning the expected performance of first-stage decisions. Neural networks are trained using recourse values computed over a sufficiently large scenario set \cite{Chou24}, enabling the model to implicitly capture stochastic variability without explicit scenario conditioning. 

During our preliminary investigations, we found that existing scenario-embedding frameworks, which have proven effective for various classical two-stage SP problems (e.g., Neur2SP \cite{Dumouchelle22}), do not extend well directly to routing problems such as the mpTSP, due to both structural and computational considerations. From a structural point of view, prior applications such as facility location problem with $n$ nodes requires only $n$ binary variables to represent first-stage decisions, whereas routing problems necessitate $\mathcal{O}(n^2)$ binary arc variables to encode both connectivity and sequence. This dimensional mismatch complicates the construction of compact, informative scenario embeddings and makes them prone to information loss. From a computational perspective, these frameworks typically rely on substantial offline data generation and training. While such costs are common and well justified in domains such as natural language processing \cite{Qin26}, where models are trained once and reused across a large number of instances with fast inference, their benefits are more context-dependent in operations research settings. In many industrial engineering applications, problem instances are solved in a more instance-specific manner, and efficient heuristics already provide strong performance. For example, the data generation and training phases of Neur2SP, on a relatively small facility location instances (e.g., 10 customers and 10 candidate locations), take several minutes \cite{Dumouchelle22} on high-performance computing resources with extensive parallelization (e.g., 43 parallel processes), whereas well-designed heuristics can produce high-quality solutions within seconds \cite{Almeida24}. As learning-based methods are heuristic in nature, they are more appropriately compared with other heuristics rather than exact algorithms. As a result, the trade-off between offline training cost and practical benefit is less straightforward. Motivated by this observation, we aim to develop an approach that leverages offline computation more efficiently, enabling practical deployment without incurring excessive training costs and risk of overfitting.

The remainder of this paper is organized as follows. Section \ref{sec:lr} reviews the relevant literature on stochastic programming and machine learning approaches in routing and combinatorial optimization. Section 3 formulates the mpTSP as a two-stage stochastic programming model and introduces a Deterministic Equivalent formulation, which serves as the baseline for comparison. Section 4 presents the proposed neural network-based surrogate modeling approaches for the recourse problem, detailing their architectures and training procedures. Section 5 describes the experimental setup, evaluation metrics, and computational results. Finally, Section 6 concludes the paper and outlines potential directions for future research.

\section{Literature Review}\label{sec:lr}
Early work on the mpTSPs focused on capturing realistic uncertainty in travel times while maintaining computational tractability. In the initial formulation presented in \cite{Tadei17}, stochastic travel times are modeled via a finite set of scenarios derived from real traffic data. Each arc has several possible travel times depending on congestion levels, and scenario probabilities are estimated empirically from observed frequencies. 
To address the resulting computational burden that is commonly encountered in stochastic routing problems (e.g., \cite{Yu12,Chou21}), a scenario reduction procedure is applied. The resulting stochastic problem is converted into a \textit{Deterministic approximation} and the mixed-integer programming model (MILP) is solved using the Concorde TSP solver \cite{Applegate06}. Built on this framework, the two-stage formulation of mpTSP is introduced in \cite{Perboli17}, where the first-stage decision is to construct a feasible tour and the second-stage recourse problem determines the paths to choose under scenario realizations. The relevance of adaptive and dynamic decision-making in stochastic routing contexts has been emphasized in \cite{Psaraftis95}. To solve the two-stage model, \cite{Perboli17} propose both subtour-elimination and flow-based formulations, and employ the Progressive Hedging algorithm \cite{Rockafellar91} to enhance scalability.

Further extensions of the stochastic TSP framework often focus on exploring alternative representations of uncertainty \cite{Toriello14, Montemanni07, Chou21, Fadda20}. A dynamic TSP with stochastic arc costs is proposed in \cite{Toriello14}, the aim is to enable sequential decision-making as travel times are realized over time. In a more recent work \cite{Fadda20}, the mpTSP is further extended by incorporating dependencies between arc costs, thereby modeling more realistic stochastic structures. These extended formulations, however, also further amplify computational complexity compared to prior modelling. All studies in the literature demonstrate that the computational burden of stochastic mpTSPs grows rapidly with the number and complexity of scenarios increase, motivating the exploration of alternative solution approaches.

Recent advances in machine learning for large-scale optimization have raised increasing interest in using predictive models to approximate complex optimization problems \cite{Abbasi20,Yang22}. Within the field of Stochastic Programming, several recent surveys and methodological contributions have explored the integration of learning-based models with stochastic programming frameworks \cite{Dumouchelle22,Chou25,Chou24,Chou23,Sadana25}. In two-stage stochastic programming, the second-stage recourse problem is itself an optimization problem nested within the first-stage problem. This structure makes direct prediction of complete solutions particularly challenging, due to the highly nonlinear and high-dimensional mapping between first-stage decisions and global optimal outcomes. Therefore, a particularly promising approach is to employ neural networks as surrogate models to approximate only the second-stage recourse problem \cite{Dumouchelle22, Chou25, Chou24}. Instead of directly approximating the expected recourse function value, a quantile-based neural network framework is proposed in \cite{Alcantara25} that learns specific quantiles of the distribution of the recourse problem. This methodology provides a more detailed statistical characterization of uncertainty compared to expectation-based approaches. Other uses of neural networks in stochastic or robust optimization can be found in the existing literature \cite{Goerigk23,Angioni25}.

Overall, recourse function approximation \cite{Chou24,Chou25} is most effective when the recourse problem is relatively easy to solve, even if the complete two-stage SP problem is computationally demanding. This property makes the approach well-suited for the mpTSP, motivating its investigation in the present study.

\section{Problem Formulation}\label{sec:pb}
\subsection{Two-stage mpTSP}
Let $V=\{1,\dots,N\}$ denote a set of nodes in a graph including the depot, and let $E=\{(i,j): i,j\in V, i\neq j\}$ denote the set of arcs between nodes. For each pair of nodes $(i,j)\in E$, a set of alternative paths $P_{ij}$ is available, representing different feasible paths (e.g., road segments or congestion patterns) between the two nodes. The travel time associated with each path is subject to stochastic variation due to uncertain traffic conditions. The objective of the mpTSP is to determine a ``here-and-now'' optimal feasible tour that minimizes the sum of fixed travel costs and the expected recourse costs arising from ``wait-and-see'' adjustments to traffic conditions captured in the second-stage problem. The general representation of a two-stage stochastic programming formulation of the mpTSP with stochastic travel costs can be represented as follows:
\begin{align}\label{eq:sp_1}
\min_{x \in \mathcal{X}} c^Tx +  \mathbb{E}_{\xi}[Q(x, \xi)]
\end{align}
where:
\begin{itemize}
    \item $x = \{x_{ij}\}_{(i,j)\in E}$ denotes first-stage decision variables, with $x_{ij} = 1$ if arc $(i,j)$ is included in the tour and $x_{ij} = 0$ otherwise;
    \item $c = \{c_{ij}\}_{(i,j)\in E}$ represents the fixed costs of the arcs in the graph;
    \item $\xi \in \Xi$ is the set of random parameters representing the uncertainty in travel times, assumed to follow certain distributions;
    \item $Q(x,\xi)$ is the \emph{recourse function}, representing the optimal second-stage cost for first-stage decisions $x$ under the uncertainty $\xi$.
\end{itemize}

In mpTSP, the feasible set $\mathcal{X}$ enforces that $x$ defines a Hamiltonian tour, i.e., each node is visited exactly once and subtours are eliminated. For a given tour $x$, the recourse problem determines the optimal path choices between visited nodes under the uncertainty $\xi$. Let $\mathcal{Y} = \{0, 1\}^{|E| \times |P|}$ be the space of binary path selection decisions. The recourse function is defined as:
\begin{align}\label{eq:sp_2}
Q(x,\xi) = 
\min_{y \in \mathcal{Y}} 
\sum_{(i,j)\in E}\sum_{p \in {P}_{ij}} 
\Delta_{ij}^p(\xi)\, y_{ij}^p(\xi),
\end{align}
s.t.
\begin{align}\label{eq:sp_3}
& \sum_{p \in P_{ij}} y_{ij}^p(\xi) = x_{ij}, \quad \forall (i,j) \in E
\end{align}
where:
\begin{itemize}
    \item $y_{ij}^p(\xi)$ is the second-stage decision variable in the feasible region $\mathcal{Y}(\mathcal{X}, \Xi)$, restricted by the feasible region of the first-stage tour $\mathcal{X}$ and the uncertainty space $\Xi$. $y_{ij}^p(\xi)$ indicates whether path \(p \in P_{ij}\) is chosen for arc \((i,j)\) under the uncertainty \(\xi\), with $y_{ij}^p(\xi)=1$ if path $p$ is selected and $y_{ij}^p(\xi)=0$ otherwise;
    \item $\Delta_{ij}^p(\xi)$ is the travel time oscillation of path $p$ on arc $(i,j)$ under the uncertainty $\xi$;
    \item Constraint (\ref{eq:sp_3}) enforces that exactly one path is selected for each chosen arc $x_{ij} = 1$, and none otherwise.
\end{itemize}

In general, problem~(\ref{eq:sp_1}) is not directly solvable. To obtain a tractable approximation, the continuous probability distribution is usually discretized into a finite scenario set $S = \{\xi^1, \dots, \xi^{|S|}\}$ with associated probabilities $\pi_s$ for each scenario $\xi^s$, satisfying \(\sum_{s \in S} \pi_s = 1\). The second-stage variables can then be defined as scenario dependent, leading to the \textit{Deterministic Equivalent} formulation of the two-stage mpTSP (see Section \ref{sec:de}). As the number of scenarios increases, the discrete approximation provides a more accurate representation of the underlying distribution $\xi$. However, problem~(\ref{eq:sp_2}) then involves $|S|$ replications of the second-stage objective $Q(x, \xi)$ and its associated feasible region $\mathcal{Y}(\mathcal{X}, \Xi)$, which can make the formulation computationally intractable when $|S|$ is large. Given the inherent uncertainty and variability in traffic conditions \cite{Susilawati13,Bauer19}, a large number of scenarios is typically required to obtain a reliable approximation. This poses a significant challenge for solving the mpTSP efficiently.

Alternative solution approaches, such as classical Benders decomposition, are not directly applicable to general two-stage stochastic mixed-integer programs, as they require the second-stage variables to be continuous \cite{Maheo24}. Although several variations of the method have been proposed to address this limitation, their application to two-stage stochastic mpTSP remains unexplored. Therefore, we adopt the Deterministic Equivalent formulation of the flow-based two-stage stochastic programming model \cite{Perboli17} as a baseline solution approach for the mpTSP, acknowledging that its computation time increases as the number of scenarios grows. 

\subsection{The Deterministic Equivalent}\label{sec:de}
The flow-based Deterministic Equivalent formulation is well-suited for mixed-integer programming solvers, as it has $\mathcal{O}(N)$ constraints, instead of the exponential number of constraints often found in sub-tour elimination based models \cite{Perboli17} .

In the flow-based model, a continuous first-stage variable $\phi_{ij}$ is introduced to represent the flow along arc $(i, j)$. For each pair of connected nodes $(i, j)$, the set of feasible paths is denoted by $P_{ij}$. Each path $p \in P_{ij}$ is characterized by a non-negative estimated mean unit travel time cost $\bar{c}_{ij}$ and a non-negative random unit travel time cost $c_{ij}^{ps}$, realized under time scenario $s \in S$. In practice, the first-stage cost $\bar{c}_{ij}$ is calculated as the Euclidean distance $d_{ij}$ between nodes divided by the expected value of the empirical velocity profile distributions observed across feasible paths. The realized travel time cost for path $p \in P_{ij}$ under scenario $s \in S$ is similarly approximated as the Euclidean distance $d_{ij}$ divided by the empirically observed velocity $v_{ij}^{ps}$ from sensor data along that path. The second-stage cost is then defined as the deviation between the realized travel time and the expected travel time for a given path $p$ and scenario $s$:
\begin{align}
\Delta_{ij}^{ps} = c_{ij}^{ps} - \bar{c}_{ij} = \frac{d_{ij}}{v_{ij}^{ps}}-\mathbb{E}_{s\in S}[\frac{d_{ij}}{\mathbb{E}_{p \in P_{ij}}[v_{ij}^{ps}]}]
\end{align}

By incorporating the scenario-dependent second-stage decision variable $y_{ij}^{ps}$, the \textit{Deterministic Equivalent} (DE) formulation of the mpTSP can be defined as follows:
\begin{align}
&\min \left[ \sum_{(i,j) \in E} \bar{c}_{ij} x_{ij} + \sum_{s \in S}\pi_s\sum_{(i,j) \in E} \sum_{p \in P_{ij}} \Delta_{ij}^{ps} y_{ij}^{ps} \right]\label{eq:DE-obj} 
\end{align}
s.t.
\begin{align}
&\sum_{\substack{j \in N \\ j \neq i}} x_{ij} = 1, \quad \forall i \in N \label{eq:DE-constraint-0}\\
&\sum_{\substack{i \in N \\ i \neq j}} x_{ij} = 1, \quad \forall j \in N \label{eq:DE-constraint-1}\\
& \sum_{\substack{i \in N \\ i \ne j}} \phi_{ij} -1 = \begin{cases} 0 , & \text{ if }j=1, \\ \sum_{\substack{k \in N \\ k \ne j}} \phi_{jk}, & \text{ otherwise} \end{cases} \quad \forall j \in N \label{eq:DE-constraint-2}\\
&\sum_{\substack{k \in N \\ k \neq 1}} \phi_{1k} = |N|\label{eq:DE-constraint-4} \\
&\phi_{ij} \leq |N| x_{ij}, \quad \forall i \in N, j \in N \label{eq:DE-constraint-5}\\
&\sum_{p \in P_{ij}} y_{ij}^{ps} = x_{ij}, \quad \forall i \in N, j \in N, s \in S \label{eq:DE-constraint-6}\\
&x_{ij}, y_{ij}^{ps} \in \{0, 1\}, \quad \forall p\in P_{ij}, i \in N, j \in N, s \in S \label{eq:DE-constraint-7}
\end{align}

The objective function (\ref{eq:DE-obj}) minimizes the total expected cost, which is composed of the deterministic travel costs from the first-stage routing and the expected stochastic travel time oscillations incurred by the adaptive second-stage path selections. Constraints (\ref{eq:DE-constraint-0})-(\ref{eq:DE-constraint-1}) enforce the fundamental TSP requirement that every node must be entered and left exactly once. Constraints (\ref{eq:DE-constraint-2})-(\ref{eq:DE-constraint-5}) ensure route feasibility by maintaining flow balance and preventing the formation of subtours. Coherence between the two stages is maintained by Constraints (\ref{eq:DE-constraint-6}), which guarantee that for every arc selected in the first-stage tour, exactly one specific path is chosen under each realized scenario. Finally, Constraints (\ref{eq:DE-constraint-7}) define the binary nature of the decision variables.

\subsection{The Deterministic Approximation}\label{sec:da}

To further motivate the choice of using a two-stage SP formulation instead of solving a deterministic TSP with approximated stochastic arc costs, in this section we provide an analysis comparing it to a \textit{Deterministic approximation} (DA) approach used in \cite{Tadei17}.

The DA method simplifies the mpTSP by substituting all stochastic parameters with their expected values (e.g., mean travel times), thereby converting the problem into a single deterministic TSP \cite{Tadei17}. Although computationally much simpler, this approach fails to account for the variability in the stochastic environment, potentially leading to solutions that perform poorly under uncertainty. We illustrate this limitation by considering a simple arc $(A, B)$ with two alternative paths and two equally probable scenarios, each associated with the travel costs:
\begin{center}
\begin{tabular}{c|cc}
\toprule
Scenario & Path 1 & Path 2 \\
\midrule
$s_1$ & 2 & 8 \\
$s_2$ & 8 & 2 \\
\bottomrule
\end{tabular}
\end{center}

Under the DA approach, the expected cost for both Path 1 and Path 2 is computed as the average, which is $5$ in both cases. Consequently, the deterministic arc cost is $c_{\textsc{DA}}(A, B) = 5$. In contrast, the two-stage SP formulation allows for adaptive second-stage decisions to select the minimum cost path upon scenario realization (as defined in Equation \ref{eq:sp_2}). In this specific instance, the SP model would always select the path with the minimum travel time ($2$) for each scenario, resulting in an expected cost of $c_{\textsc{DE}}(A, B) = 2$. This divergence in expected arc cost estimation demonstrates how the two approaches can lead to significantly different optimal TSP tours.

\begin{figure}[t]
    \centering
    \begin{subfigure}{0.47\textwidth}
        \centering
        \includegraphics[width=\textwidth]{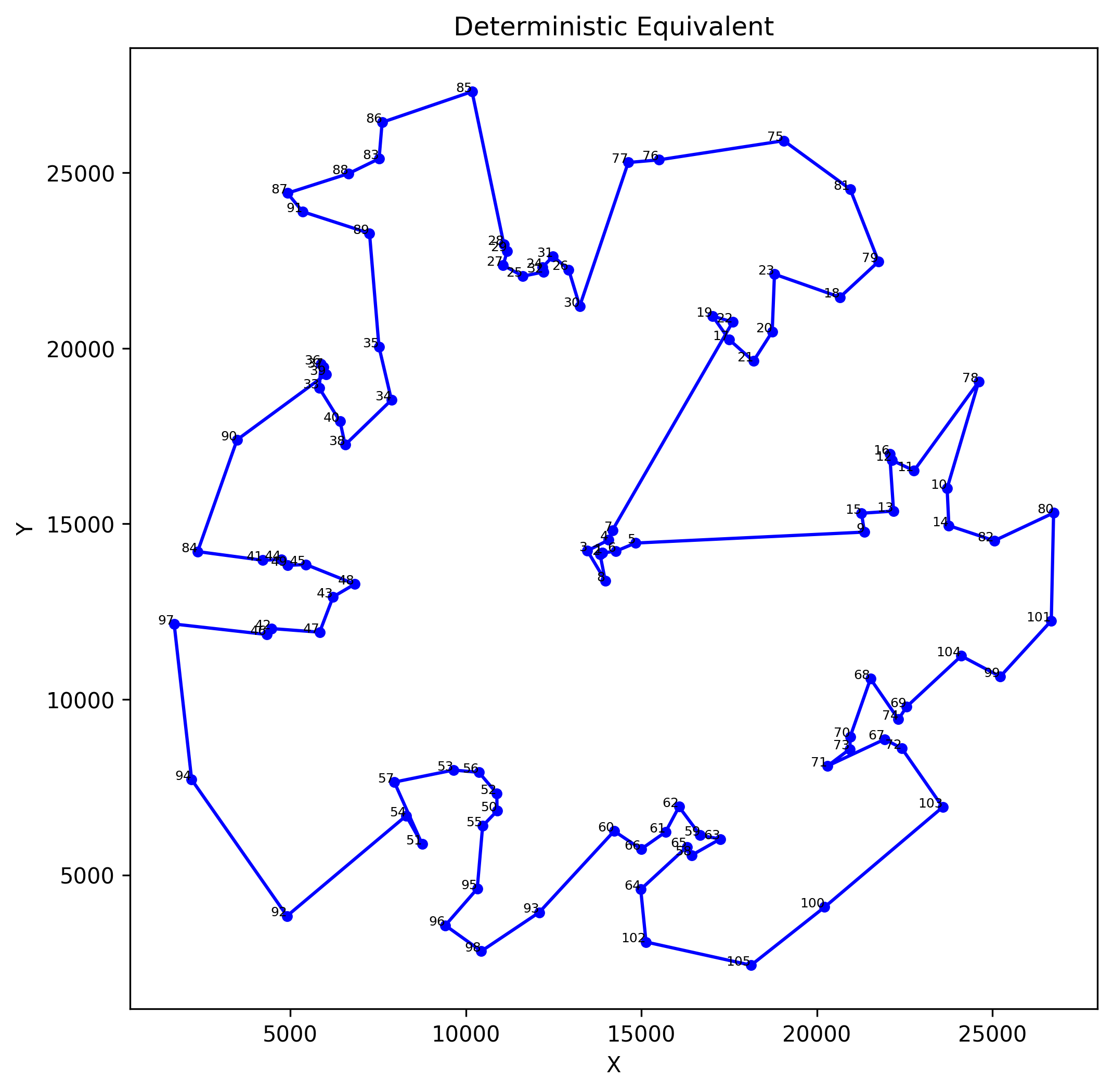}
        \caption{Optimal TSP route found by DE.}
        \label{fig:DE_route}
    \end{subfigure}
    \hfill
    \begin{subfigure}{0.47\textwidth}
        \centering
        \includegraphics[width=\textwidth]{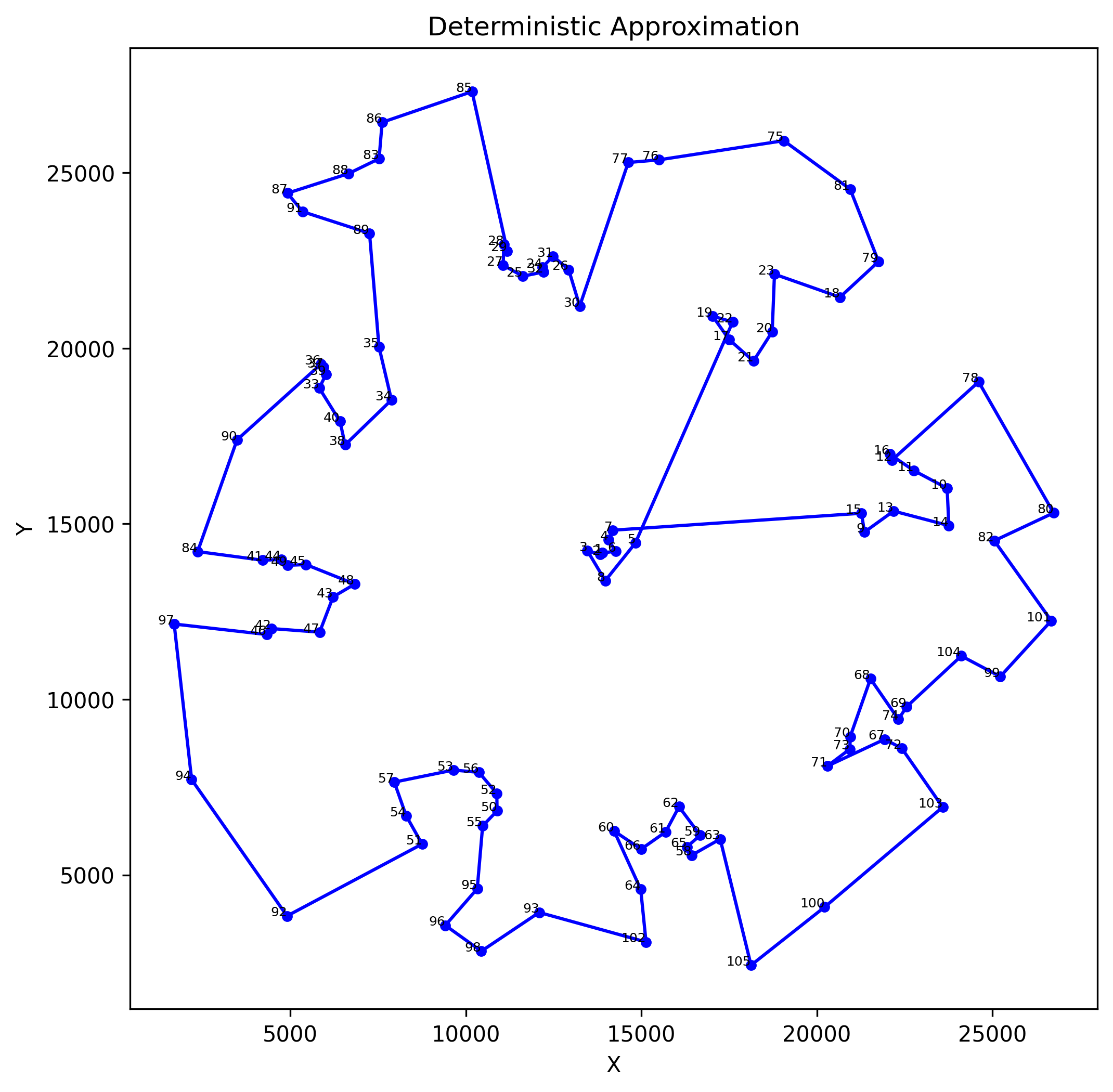}
        \caption{Optimal TSP route found by DA.}
        \label{fig:DA_route}
    \end{subfigure}
    \caption{Comparison of optimal TSP routes for a real-world instance obtained using the \textit{Deterministic Equivalent} (DE) and \textit{Deterministic Approximation} (DA) approaches}
    \label{fig:DE_DA_routes}
\end{figure}

The DA approach could provide reasonable results if travel times are symmetrically distributed with low variance (e.g., under a normal distribution). However, real-world travel times frequently exhibit asymmetric or heavy-tailed distributions \cite{Bauer19,Higatani09}, making the mean value a poor representative of true variability. Consequently, the SP formulation produces solutions that are more flexible and risk-aware, particularly when uncertainty is substantial or when paths display asymmetric uncertainty profiles. This difference is empirically validated on a real-world problem instance in Figure \ref{fig:DE_DA_routes}. The optimal DA route, obtained by solving the approximated deterministic TSP using the Lin-Kernighan-Helsgaun (LKH) heuristic \cite{LKH}, is compared against the optimal DE solution derived from the SP formulation described in Section \ref{sec:de}. When evaluated against the same set of 100 scenarios, the DE route achieves a lower total expected cost than the DA route, demonstrating the benefits of explicitly modeling uncertainty in the mpTSP.

\section{Neural Network surrogate for the mpTSP}\label{sec:nn}
To address the computational challenges arising from a large number of scenarios, we employ neural network models to approximate the recourse problem of the mpTSP, using them as surrogate models within the original optimization framework.

\subsection{General framework}
The fundamental objective is to train neural network models to approximate and replace the expected value of the recourse problem $\bar{Q}(x) = \mathbb{E}_\xi[Q(x,\xi)]$. The resulting optimization problem of (\ref{eq:sp_1}) can be expressed as:
\begin{align}\label{eq:sp_surrogate}
    \min \sum_{i=1}^n \sum_{j=1}^n\bar{c}_{ij}x_{ij} + \hat{y}
\end{align}
where $\hat{y}$ denotes the output of a trained neural network $\hat{Q}(x; \theta)$, which serves as a data-driven surrogate for the expected second-stage cost.

\begin{figure}[ht]
    \centering
    \includegraphics[width=0.8\textwidth]{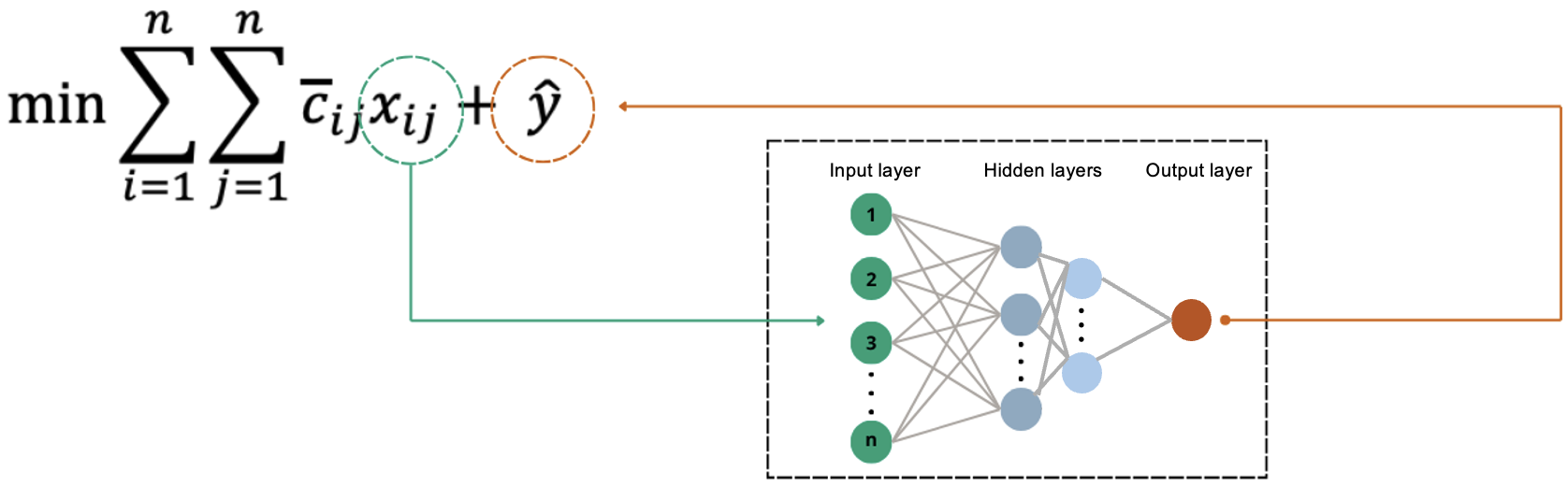}
    \caption{General framework of the neural network surrogate model in the mpTSP}
    \label{fig:nn_structure}
\end{figure}

By employing the surrogate model, the explicit recourse formulation and constraints in (\ref{eq:sp_2})-(\ref{eq:sp_3}) are replaced by an implicit approximation learned by the neural network. Instead of explicitly solving the second-stage subproblems, the network captures the relationship between first-stage decision variables and their corresponding second-stage costs, enabling efficient evaluation within the optimization model. Consequently, the overall model performance depends critically on the predictive accuracy and generalization capability of the surrogate model. Figure~\ref{fig:nn_structure} illustrates the general framework, and the following section details its design for the mpTSP.

\subsubsection{Route encoding}
\begin{figure}[t]
    \centering
    \includegraphics[width=0.55\textwidth]{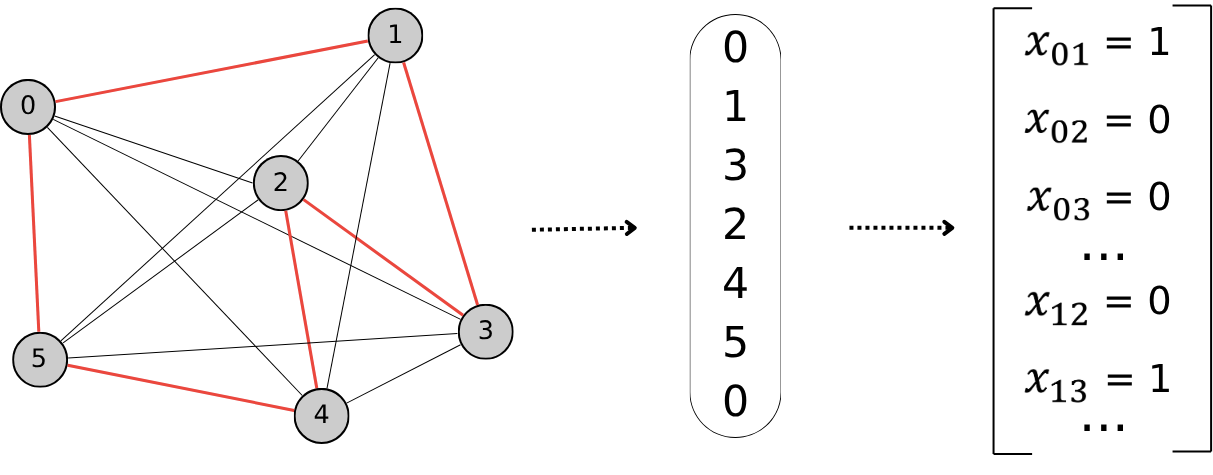}
    \caption{Illustration of the route encoding used as the input layer of the neural network}
    \label{fig:route_encoding}
\end{figure}

The surrogate model requires a compact and informative representation of combinatorial routes as input. We adopt a binary edge-incidence vector \(x\in\{0,1\}^{|E|}\), where each element indicates whether the corresponding edge is included in the route (see Figure \ref{fig:route_encoding}). This encoding effectively preserves the network's topological structure while remaining independent of the underlying neural architecture. While there exist alternative representations such as permutation vectors or adjacency matrices, the binary edge-incidence format integrates naturally with fully connected neural networks and aligns well with linear and integer programming formulations used in the underlying optimization model (see Section \ref{sec:integration}).

\subsubsection{Network architecture and training}
The training dataset $\mathcal{D}=\{(x^k,\bar{Q}(x^k))\}_{i=1}^{K}$ is constructed by sampling $K$ random route configurations $x^{k} = \{x_{ij}\}_{(i,j)\in E}$ and evaluating their corresponding expected recourse costs \(\bar{Q}(x^{k})\) over a predefined scenario set \(S = \{\xi^{1}, \ldots, \xi^{|S|}\}\). Each pair $(x^k,\bar{Q}(x^k))$ serves as an input-output example (illustrated by the green and red points in Figure~\ref{fig:nn_structure}) for training the neural network $\hat{Q}(x; \theta)$. Feed-forward neural networks (FFNN) with multiple hidden layers and rectified linear unit (ReLU) activation functions are employed as the baseline surrogate architecture in this work. This simple structure provides sufficient approximation capability for the instances investigated. The forward propagation through the network is defined as:
\begin{align}
    z_l &= w_l * \sigma(z_{l-1}) + b_l
\end{align}
where $z_l$ denotes the linear transformation at layer $l$,  $\sigma$ is the activation function, $w_l$ and $b_l$ represent the layer's weight matrix and bias vector, respectively. In this case, the activation function is given by $\text{ReLU}(a)=\max\{0,a\}$, for $a \in \mathbb{R}$. The neural network parameters $\theta = \{w_l, b_l\}$ are then learned by minimizing the mean squared error (MSE), mean absolute error (MAE) or mean absolute percentage error (MAPE), depending on the model configuration:
\begin{align}
    \min_{\theta} \; \frac{1}{K} \sum_{k=1}^{K} (\bar{Q}(x^{k}) - \widehat{Q}(x^{k};\theta))^{2}
    \label{eq:mse}
\end{align}
\begin{align}
    \min_{\theta} \; \frac{1}{K} \sum_{k=1}^{K} |\bar{Q}(x^{k}) - \widehat{Q}(x^{k};\theta)|
    \label{eq:mae}
\end{align}
\begin{align}
\min_{\theta} \frac{100}{K} \sum_{k=1}^{K}
\left| \frac{\bar{Q}(x^{k}) - \widehat{Q}(x^{k};\theta)}{\bar{Q}(x^{k})} \right|
\label{eq:mape}
\end{align}

Once trained, the surrogate model provides an efficient approximation of the expected recourse cost, enabling fast evaluation of candidate routes and avoiding repeated scenario-based computations. In the stochastic mpTSP, where travel times are random and potentially correlated across edges, this data-driven approximation is anticipated to effectively capture structural dependencies in the recourse function and substantially reduce computational complexity of the original optimization model.

\subsubsection{Integration with the optimization model}
\label{sec:integration}
A FFNN with ReLU activation function can be embedded into a MILP by explicitly reformulating each neuron using auxiliary continuous and binary variables \cite{Fischetti2018}. For each neuron $h^\ell_j$ at position $j$ in layer $\ell$, its actual value can be represented as:
\begin{align}
h_j^\ell = \sigma\left(\sum_{i=1}^{d_{\ell-1}} w^{\ell-1}_{ij} \, h^{\ell-1}_i + b^{\ell-1}_j \right),
\end{align}
where $d_{\ell-1}$ is the dimension of the previous layer, $w^{\ell-1}_{ij}$ and $b^{\ell-1}_j$ are the pretrained weights and bias, respectively, $\sigma(\cdot)$ is the ReLU function.

To linearize the ReLU function, we introduce nonnegative variables $\hat{h}_j^\ell,\,\check{h}_j^\ell \ge 0$ and a binary activation indicator $z_j^\ell \in \{0,1\}$, and impose:
\begin{align}
\sum_{i=1}^{d_{\ell-1}} w^{\ell-1}_{ij} \, \hat{h}^{\ell-1}_i + b^{\ell-1}_j&= \hat{h}_j^\ell - \check{h}_j^\ell \label{eq:neuron1}\\
z_j^\ell = 1 &\Rightarrow \hat{h}_j^\ell = 0 \label{eq:relu1}\\
z_j^\ell = 0 &\Rightarrow \check{h}_j^\ell = 0 \label{eq:relu2}
\end{align}

This set of constraints reproduces the piecewise-linear behavior of the ReLU activation: when the value on the left-hand side of equation~(\ref{eq:neuron1}) is strictly positive, the model forces $z_j^\ell = 0$ and $\hat{h}_j^\ell > 0$; when the left-hand side is non-positive, the model forces  $z_j^\ell = 1$ and $\hat{h}_j^\ell = 0$. This approach ensures an exact MILP representation of the neural network, allowing it to be fully integrated into the original optimization problem.

The number of additional constraints and variables introduced is directly proportional to the size of the neural network. Consequently, embedding a large network can substantially increase the overall problem complexity and solution time. In practice, a trade-off must be found between the predictive accuracy of the surrogate model and the resulting computational tractability of the combined MILP.

\subsection{Choice of training paradigm}
As discussed in Section \ref{sec:intro}, two primary paradigms exist for training neural network surrogates of the recourse function in two-stage SP problems: \emph{scenario-conditioned} and \emph{decision-only}.

We first considered the \emph{scenario-conditioned} surrogate modeling paradigm following the Neur2SP framework \cite{Dumouchelle22}. A neural network predictor $\phi(x,\zeta;\theta)$ is trained to estimate the recourse function value associated with a first-stage decision $x$, conditioned on a set of scenarios $\{\xi_s\}_{s=1}^{|S|}$. With $p_s$ being the probability of scenario $\xi^s$, the expected value of the recourse function is approximated empirically over the scenario set as:
\begin{equation}
\bar{Q}(x,\{\xi_s\}_{s=1}^{|S|}) = \sum_{s=1}^{|S|}p_sQ(x, \xi_s)
\end{equation}

Following the approach in \cite{Dumouchelle22}, two distinct architectures are considered. The \textit{single-cut} architecture learns a direct mapping $(x,\{\xi^s\}_{s=1}^{|S|}) \mapsto \bar{Q}(x,\{\xi^s\}_{s=1}^{|S|})$. The entire scenario set is first encoded and aggregated into a joint representation $\zeta$, which is then concatenated with the first-stage decision vector $x$ as input to a single neural network. Alternatively, the \textit{multi-cut} architecture learns the mapping $(x,\xi^s) \mapsto Q(x,\xi^s)$ separately for each scenario $\xi^s$. A distinct neural network is trained for each scenario, and the resulting outputs are averaged to approximate the expected recourse function value. Both frameworks have proven effective in various classical two-stage SP problems, including the capacitated facility location problem and the pooling problem \cite{Dumouchelle22}. A commonality in these successful applications is that the first-stage decision space is relatively low-dimensional and primarily composed of binary variables.

In contrast, the mpTSP involves a substantially higher-dimensional decision space. While a facility location problem with $n$ nodes requires only $n$ binary variables to represent first-stage decisions, routing problems necessitate $\mathcal{O}(n^2)$ binary arc variables to encode both connectivity and order (see Figure \ref{fig:route_encoding}). This large dimensional mismatch between the first-stage solution vector $x$ and the aggregated scenario embeddings $\zeta$ makes their direct concatenation as input impractical for the mpTSP, often leading to poor network training performance. Furthermore, travel time uncertainty in routing problems is both high-dimensional and frequently correlated across the network. This makes it challenging to construct compact and informative scenario embeddings without significant information loss. Minor perturbations in the first-stage solution can also induce large variations in the recourse cost due to topology-dependent route changes. For these reasons, in this work we adopt an alternative surrogate modeling approach that does not rely on scenario embeddings.

The \textit{decision-only} surrogate modeling paradigm works under the assumption that, given a sufficiently large scenario set, the empirical estimation of the recourse function value converges to its true expectation and thus provides a stable performance indicator for any first-stage solution \cite{Chou24}. Formally, the expected value of the recourse function can be approximated as 
\begin{align}
 \bar{Q}(x)\approx\lim_{|S| \rightarrow +\infty }\mathbb{E}_{s \in S} [Q(x,\xi_s)] 
\end{align}

Formalized in the NN2SP framework \cite{Chou25,Chou24}, this approach aims to obtain a compact representation that captures the expected performance of first-stage solutions. This avoids the computationally expensive process of repeated scenario-based evaluation within the SP optimization model. The NN2SP framework has demonstrated efficacy in the single-source facility location problem \cite{Chou24} and has been preliminarily evaluated on the mpTSP \cite{Chou25}. In the present study, we fully adopt this framework for the mpTSP by training a surrogate model that directly maps each feasible tour to its expected recourse cost, i.e., $x \mapsto \bar{Q}(x)$.

In the mpTSP, solving the recourse problem is straightforward, as it involves selecting the shortest path among the available options for the chosen edges. To avoid performing this computation repeatedly within the optimization model, it is approximated by a neural network that serves as a surrogate model. The training dataset consists of pairs $(x, \bar{Q}(x))$, where each $\bar{Q}(x)$ is computed over a large scenario set to ensure convergence. This convergence is empirically verified, confirming that the estimated values remain stable as additional scenarios are included. Since the recourse problems are independent across scenarios and first-stage solutions, the training dataset can be efficiently generated in parallel to reduce computational time.

Once trained, the neural network $\widehat{Q}(x;\theta)$ provides a compact representation of the estimated performance of first-stage solutions and serves as a surrogate model for $\mathbb{E}_{\xi}[Q(x,\xi)]$ in the original optimization problem. This substitution provides a tractable and robust approximation of the expected recourse function value, enabling efficient stochastic optimization for large-scale mpTSP instances.

\section{Experimental Results}\label{sec:exp}
\subsection{Dataset and settings}\label{sec:data}

The proposed method is evaluated on two different problem instances drawn from the mpTSP benchmark datasets previously used in \cite{Tadei17, Perboli17}. Both instances are defined with 100 discrete scenarios, representing all possible values for the stochastic travel time variable. For each pair of nodes, three alternative paths are available. Additional instance details are as follows:

\begin{itemize}
    \item \texttt{MPTSPs\_D3\_50}: 56 customers distributed evenly across city center and suburban area.
    \item \texttt{MPTSPs\_D1\_100}: 105 customers situated exclusively in suburban area.
\end{itemize}

Training datasets $\mathcal{D}=\{(x^k,\bar{Q}(x^k))\}_{i=1}^{K}$ were constructed by randomly sampling first-stage solutions $x^k$ and computing the corresponding expected recourse function value $\bar{Q}(x^k)$ over a scenario set $S$. Given that the total number of possible tours in an undirected TSP with $n$ nodes is $\frac{n!}{2}$, the training set cannot provide an exhaustive coverage of the solution space. To assess the impact of training set size on surrogate model performance, we generated $K$=10,000 and $K$=50,000 training observations for the 56-node instance, and $K$=100,000 and $K$=200,000 training observations for the 105-node instance. For comparison with the baseline DE approach, we constructed scenario sets with cardinalities $|S| \in \{3, 10, 20, 40, 60, 80, 100\}$.

To account for stochastic variability and ensure robust results, for all the experiments, we executed 30 independent runs and calculated the average performance metrics across all runs. Accordingly, 30 statistically independent training datasets were generated for each unique combination of problem size $n=\{56,105\}$, training set size $K$, and scenario set size $|S|$. 

To guarantee the independency of the different runs, the original 100 scenarios given in the benchmark dataset \cite{mpTSP_benchmark} were expanded to 1,000 independent scenarios sampled based on the empirical distribution of the observed travel times. This 1,000-scenario set was partitioned into two disjoint pools $\{\mathcal{S}_{\text{T}},\mathcal{S}_{\text{OOS}}\}$ of 500 scenarios each, designated for training and out-of-sample testing, respectively. Identical scenario sets were used for both the NN surrogate model and the DE approaches across all experiments to ensure a fair and consistent comparative basis.

Neural networks were trained using the Adam optimizer with an initial learning rate $\eta = 0.01$, minimizing the Mean Absolute Error as loss function. To prevent overfitting, $\mathtt{EarlyStopping}$ with a patience of 5 epochs was employed to terminate the training process upon observing a lack of improvement in the validation loss. Additionally, the $\mathtt{ReduceLROnPlateau}$ callback was applied to decrease the learning rate by a factor of $0.5$ (minimum rate $10^{-6}$) when the validation loss plateaued. The training proceeded for a maximum limit of 1,000 epochs.

All the preliminary experiments reported in this work were executed on a MacBook Air equipped with an Apple M4 chip and 16 GB of RAM. The optimization problems were implemented using Python 3.10.6 and solved with the Gurobi Optimizer 12.0.1. The neural networks were developed and trained using Keras 3.9.2 with a TensorFlow backend (version 2.19.0). All computations were performed locally without GPU acceleration or cloud-based resources.

\subsection{Parameter tuning}

For each combination of problem size $n$ and training set size $K$, a grid search was performed to identify the most suitable neural network architecture for the task. The search explored network configurations with two or three hidden layers, each containing 8 or 16 neurons, and batch sizes of 256, 512, and 1024. Larger architectures were excluded due to the computational complexities discussed in Section~\ref{sec:integration}. The mean absolute error was monitored for \texttt{EarlyStopping}, and the candidate achieving the lowest expected total cost when embedded as a surrogate model in the mpTSP and evaluated on a 500 scenario out-of-sample test set ($\mathcal{S}_{\text{OOS}}$) was selected as the optimal architecture. The resulting optimal configurations were then used to train all surrogate models for the corresponding problem instances, as summarized in Table~\ref{tab:selected_architectures}.

\begin{table*}[t]
\centering
\caption{Selected Neural Network Architectures Following Parameter Tuning}
\label{tab:selected_architectures}
\begin{tabular}{cccc}
\toprule
\textbf{Nodes} & \textbf{Observations} & \textbf{Layers} & \textbf{Batch Size} \\
\midrule
56   & 10,000  & [16, 16]    & 256 \\
56   & 50,000  & [8, 8, 16]  & 1024 \\
105  & 100,000 & [8, 16]     & 1024 \\
105  & 200,000 & [8, 8, 16]  & 1024 \\
\bottomrule
\end{tabular}
\end{table*}

\subsection{Predictive performance of the neural networks}

The predictive accuracy of the trained models was quantified using standard regression metrics: MAE, MAPE, and $R^2$. These metrics were evaluated on both the training and out-of-sample test sets. Additionally, we report the data generation time and the neural network training time. Table \ref{tab:nn_metrics_all} presents the average performance and execution times across 30 independent runs for each scenario set size.

From Table~\ref{tab:nn_metrics_all}, the training performance exhibits high accuracy across all instances. In particular, the $R^2$ values are consistently close to 1.0, indicating an excellent fit and an effective capture of the variability in the recourse function cost. The $\text{MAPE}$ of the out-of-sample test set remains below 0.1\%, demonstrating highly accurate approximations of the expected cost. The consistency of performance metrics across both training and test sets, as well as across different scenario sizes, suggests minimal overfitting and strong generalization capability to unseen data.

Given that the number of observations is small relative to the size of the solution space, predictive accuracy improves as the training set increases. Regarding computational performance, data generation time scales linearly with the number of observations. Under the current settings, where the mpTSP recourse problem is relatively easy to solve, the average data generation time per observation is approximately $2\times10^{-4}$ seconds for instances with 56 nodes and $6\times10^{-4}$ seconds for instances with 105 nodes, with limited sensitivity to the size of the scenario set within the tested range. Training time increases with both problem size and dataset scale, ranging from a few seconds for small instances to approximately two minutes for larger datasets, while remaining computationally manageable.

Overall, these results confirm that the proposed small neural networks can reliably approximate the mapping $x \mapsto \mathbb{E}[Q(x,\xi)]$ of the recourse function value, making them suitable for subsequent optimization tasks as surrogate models.

\begin{table*}[!t]
\centering
\caption{Neural Network Training Performance Metrics (Averaged Over 30 Runs)}
\label{tab:nn_metrics_all}

\resizebox{\textwidth}{!}{
\begin{tabular}{c|ccc|ccc|cc}
\multicolumn{9}{c}{\textbf{Problem instance: 56 nodes - 10,000 observations}\vspace{1mm}} \\
\toprule
\textbf{Scenario} & \multicolumn{3}{c|}{\textbf{Train}} & \multicolumn{3}{c|}{\textbf{Test}} & \textbf{Data Gen.} & \textbf{NN Train} \\
Set & MAE & MAPE (\%) & $R^2$ & MAE & MAPE (\%) & $R^2$ & \textbf{Time (s)} & \textbf{Time (s)} \\
\midrule
3   & 83.20 & 0.10 & 0.9996 & 82.31 & 0.10 & 0.9996 & 1.92 & 3.43 \\
10  & 74.95 & 0.09 & 0.9995 & 74.98 & 0.09 & 0.9995 & 1.95 & 3.25 \\
20  & 81.54 & 0.10 & 0.9994 & 83.09 & 0.10 & 0.9993 & 1.97 & 2.96 \\
40  & 79.11 & 0.10 & 0.9994 & 79.70 & 0.10 & 0.9994 & 2.05 & 3.06 \\
60  & 75.96 & 0.09 & 0.9995 & 74.84 & 0.09 & 0.9995 & 2.10 & 3.23 \\
80  & 74.93 & 0.09 & 0.9995 & 74.51 & 0.09 & 0.9995 & 2.18 & 3.64 \\
100 & 78.61 & 0.10 & 0.9994 & 80.53 & 0.10 & 0.9993 & 2.25 & 3.46 \\
\bottomrule
\end{tabular}}
\vspace{0.5em}

\resizebox{\textwidth}{!}{
\begin{tabular}{c|ccc|ccc|cc}
\multicolumn{9}{c}{\textbf{Problem instance: 56 nodes - 50,000 observations}\vspace{1mm}} \\
\toprule
\textbf{Scenario} & \multicolumn{3}{c|}{\textbf{Train}} & \multicolumn{3}{c|}{\textbf{Test}} & \textbf{Data Gen.} & \textbf{NN Train} \\
Set & MAE & MAPE (\%) & $R^2$ & MAE & MAPE (\%) & $R^2$ & \textbf{Time (s)} & \textbf{Time (s)} \\
\midrule
3   & 48.21 & 0.06 & 0.9998 & 48.17 & 0.06 & 0.9998 & 9.47 & 4.45 \\
10  & 57.85 & 0.07 & 0.9997 & 57.92 & 0.07 & 0.9997 & 9.49 & 5.01 \\
20  & 56.50 & 0.07 & 0.9997 & 56.59 & 0.07 & 0.9997 & 9.52 & 5.41 \\
40  & 58.53 & 0.07 & 0.9997 & 59.13 & 0.07 & 0.9997 & 9.60 & 4.74 \\
60  & 59.40 & 0.08 & 0.9996 & 60.15 & 0.08 & 0.9996 & 9.65 & 5.16 \\
80  & 60.00 & 0.08 & 0.9997 & 60.35 & 0.08 & 0.9997 & 9.75 & 3.60 \\
100 & 54.62 & 0.07 & 0.9997 & 54.40 & 0.07 & 0.9997 & 9.85 & 3.66 \\
\bottomrule
\end{tabular}
}
\vspace{0.5em}

\resizebox{\textwidth}{!}{
\begin{tabular}{c|ccc|ccc|cc}
\multicolumn{9}{c}{\textbf{Problem instance: 105 nodes - 100,000 observations}\vspace{1mm}} \\
\toprule
\textbf{Scenario} & \multicolumn{3}{c|}{\textbf{Train}} & \multicolumn{3}{c|}{\textbf{Test}} & \textbf{Data Gen.} & \textbf{NN Train} \\
Set & MAE & MAPE (\%) & $R^2$ & MAE & MAPE (\%) & $R^2$ & \textbf{Time (s)} & \textbf{Time (s)} \\
\midrule
3   & 110.31 & 0.03 & 0.9999 & 111.45 & 0.03 & 0.9999 & 62.50 & 97.57 \\
10  & 109.17 & 0.03 & 0.9999 & 110.34 & 0.03 & 0.9999 & 62.60 & 91.77 \\
20  & 110.31 & 0.03 & 0.9999 & 111.62 & 0.03 & 0.9999 & 62.70 & 98.99 \\
40  & 105.31 & 0.03 & 0.9999 & 106.21 & 0.03 & 0.9999 & 62.90 & 105.39 \\
60  & 101.56 & 0.03 & 0.9999 & 102.11 & 0.03 & 0.9999 & 63.50 & 94.04 \\
80  & 106.20 & 0.03 & 0.9999 & 107.50 & 0.03 & 0.9999 & 67.50 & 90.16 \\
100 & 102.71 & 0.03 & 0.9999 & 103.63 & 0.03 & 0.9999 & 68.00 & 95.77 \\
\bottomrule
\end{tabular}
}
\vspace{0.5em}

\resizebox{\textwidth}{!}{
\begin{tabular}{c|ccc|ccc|cc}
\multicolumn{9}{c}{\textbf{Problem instance: 105 nodes - 200,000 observations}\vspace{1mm}} \\
\toprule
\textbf{Scenario} & \multicolumn{3}{c|}{\textbf{Train}} & \multicolumn{3}{c|}{\textbf{Test}} & \textbf{Data Gen.} & \textbf{NN Train} \\
Set & MAE & MAPE (\%) & $R^2$ & MAE & MAPE (\%) & $R^2$ & \textbf{Time (s)} & \textbf{Time (s)} \\
\midrule
3   & 85.23 & 0.03 & 0.9999 & 86.38 & 0.03 & 0.9999 & 124.50 & 112.40 \\
10  & 82.47 & 0.02 & 0.9999 & 83.56 & 0.02 & 0.9999 & 124.60 & 113.70 \\
20  & 83.91 & 0.02 & 0.9999 & 84.99 & 0.02 & 0.9999 & 125.00 & 115.30 \\
40  & 80.78 & 0.02 & 0.9999 & 81.77 & 0.02 & 0.9999 & 125.20 & 118.90 \\
60  & 78.32 & 0.02 & 0.9999 & 79.33 & 0.02 & 0.9999 & 125.30 & 120.60 \\
80  & 79.45 & 0.02 & 0.9999 & 80.46 & 0.02 & 0.9999 & 125.50 & 122.80 \\
100 & 81.12 & 0.02 & 0.9999 & 82.26 & 0.02 & 0.9999 & 125.70 & 124.10 \\
\bottomrule
\end{tabular}
}
\end{table*}

\subsection{Optimization performance of the surrogate models}

The trained neural network models were subsequently integrated as surrogate models within the two-stage stochastic optimization framework. The integration was facilitated by the \texttt{Gurobi Machine Learning} Python package \cite{gurobi}, which embeds the neural network directly into the MILP model as a predictive constraint, following the MILP embedding formulation described in Section \ref{sec:integration}.

The optimal solutions obtained using the surrogate models (Section~\ref{sec:nn}) are then compared with those derived from the \emph{Deterministic Equivalent} formulation (Section~\ref{sec:de}). The comparison focuses on two primary aspects: (i) solution accuracy, evaluated on an out-of-sample dataset, and (ii) computational time. 

The cost of a given solution $x$ is evaluated over an out-of-sample scenario set as:
\[
\text{Cost}_{\text{OOS}}(x) = c^Tx+\frac{1}{|\mathcal{S}_{\text{OOS}}|} \sum_{s \in \mathcal{S}_{\text{OOS}}} Q(x, \xi_s),
\]
where $\mathcal{S}_{\text{OOS}}$ denotes a set of 500 previously unseen scenarios (Section~\ref{sec:data}), and $Q(x, \xi_s)$ represents the recourse function cost of solution $x$ under scenario $s$.

Performance differences are quantified using the relative cost gap between the optimal solution $x^*$ found by the surrogate model (NN) and the benchmark solution from the Deterministic Equivalent (DE) model:
\[
\text{GAP} (\%) = \frac{\text{Cost}_{\text{OOS}}(x^*_\text{NN}) - \text{Cost}_{\text{OOS}}(x^*_\text{DE})}{\text{Cost}_{\text{OOS}}(x^*_\text{DE})} \times 100\%
\]
where negative values indicate that the surrogate model achieves a lower total cost than the DE benchmark.

Table \ref{tab:oos_all_combined} report the quality of the optimal solutions ($x^*$), obtained using both the neural network surrogate and the deterministic equivalent model, evaluated over the out-of-sample scenario set $\mathcal{S}_{\text{OOS}}$. For each scenario set size, the reported values are the average cost $\text{Cost}_{\text{OOS}}(x^*)$ and the corresponding relative cost gap $\text{GAP} (\%)$ computed over 30 independent experimental runs.

From the table we can observe that, across all configurations, the NN surrogate consistently identifies optimal solutions with out-of-sample costs comparable to, and in some cases better than, those obtained with the DE model. The small relative gaps indicate that the surrogate is capable of to producing near-optimal solutions consistently via learned generalization.

With respect to training dataset size, strong performance is maintained across both smaller and larger datasets, suggesting that solution quality does not scale linearly with the number of training samples. In particular, models trained on smaller datasets (10,000 or 100,000 observations) exhibit marginally better average performance, as reflected by slightly lower out-of-sample costs and a higher frequency of negative cost gaps.

This behavior provides further insight into the role of the NN surrogate. Rather than simply approximating first-stage solutions and their corresponding recourse values, the model appears to capture underlying structural relationships between decision variables, enabling a meaningful degree of generalization. As a result, even when trained on relatively limited data, the surrogate remains capable of effectively distinguishing between high and low quality first-stage solutions. In contrast, the performance of the DE approach is strictly tied to the selected scenario set, which may limit its robustness with respect to unseen realizations.

\begin{table*}[ht]
\centering
\caption{Out-of-Sample Evaluation of Optimal Solutions ($x^*$) Obtained using the NN Surrogate and the DE Model
}
\label{tab:oos_all_combined}
\caption*{\texttt{Problem instance: MPTSPs\_D3\_50 (56 nodes)}}

\resizebox{\textwidth}{!}{
\begin{tabular}{c| ccc | ccc}
\toprule
& \multicolumn{3}{c|}{{$K$=10,000}} & \multicolumn{3}{c}{{$K$=50,000}} \\
\textbf{Scenario} & \textbf{DE} & \textbf{NN Sur.} & \textbf{GAP (\%)} & \textbf{DE} & \textbf{NN Sur.} & \textbf{GAP (\%)} \\
\midrule
3 & 17219.60 & 16417.00 & -4.66 & 17144.39 & 16726.77 & -2.44 \\
10 & 15534.21 & 15288.04 & -1.58 & 15517.79 & 15447.64 & -0.45 \\
20 & 15238.17 & 15147.95 & -0.59 & 15226.44 & 15215.76 & -0.07 \\
40 & 14955.52 & 14950.56 & -0.03 & 14925.12 & 14974.75 & +0.33 \\
60 & 14950.12 & 14899.77 & -0.34 & 14909.49 & 14926.92 & +0.12 \\
80 & 14882.52 & 14867.68 & -0.10 & 14840.60 & 14871.34 & +0.21 \\
100 & 14854.10 & 14887.10 & +0.22 & 14867.71 & 14860.93 & -0.05 \\
\bottomrule
\end{tabular}
}

\caption*{\texttt{Problem instance: MPTSPs\_D1\_100 (105 nodes)} }
\resizebox{\textwidth}{!}{
\begin{tabular}{c | ccc | ccc}
\toprule
& \multicolumn{3}{c|}{$K$=100,000} & \multicolumn{3}{c}{$K$=200,000} \\
\textbf{Scenario} & \textbf{DE} & \textbf{NN Sur.} & \textbf{GAP (\%)} & \textbf{DE} & \textbf{NN Sur.} & \textbf{GAP (\%)} \\
\midrule
3 & 26350.17 & 25193.08 & -4.39 & 26513.93 & 25526.50 & -3.72 \\
10 & 24389.01 & 24340.48 & -0.20 & 24327.76 & 24022.63 & -1.25 \\
20 & 23829.08 & 23703.11 & -0.53 & 23825.58 & 23804.62 & -0.09 \\
40 & 23508.93 & 23475.38 & -0.14 & 23664.45 & 23592.40 & -0.30 \\
60 & 23464.66 & 23460.36 & -0.02 & 23440.02 & 23504.10 & +0.27 \\
80 & 23415.00 & 23436.47 & +0.09 & 23410.13 & 23505.80 & +0.41 \\
100 & 23367.51 & 23412.37 & +0.19 & 23373.21 & 23469.67 & +0.41 \\
\bottomrule
\end{tabular}
}
\end{table*}

\begin{figure}[t]
    \centering
    \begin{subfigure}[b]{0.49\textwidth}
        \centering
        \includegraphics[width=\textwidth]{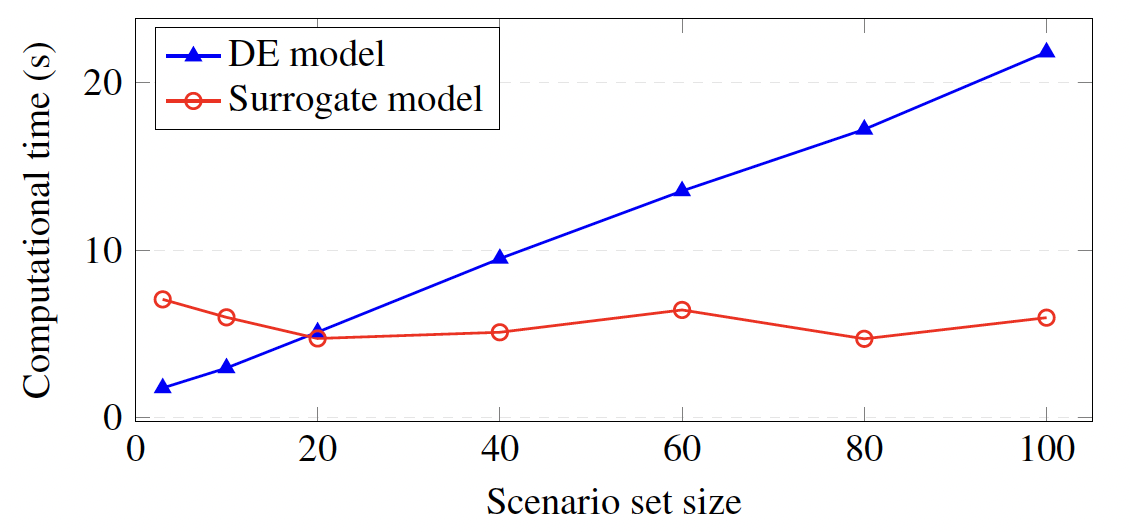}
        \caption{56 nodes - 10,000 observations}
        \label{fig:sub1}
    \end{subfigure}
    \hfill
    \begin{subfigure}[b]{0.49\textwidth}
        \centering
        \includegraphics[width=\textwidth]{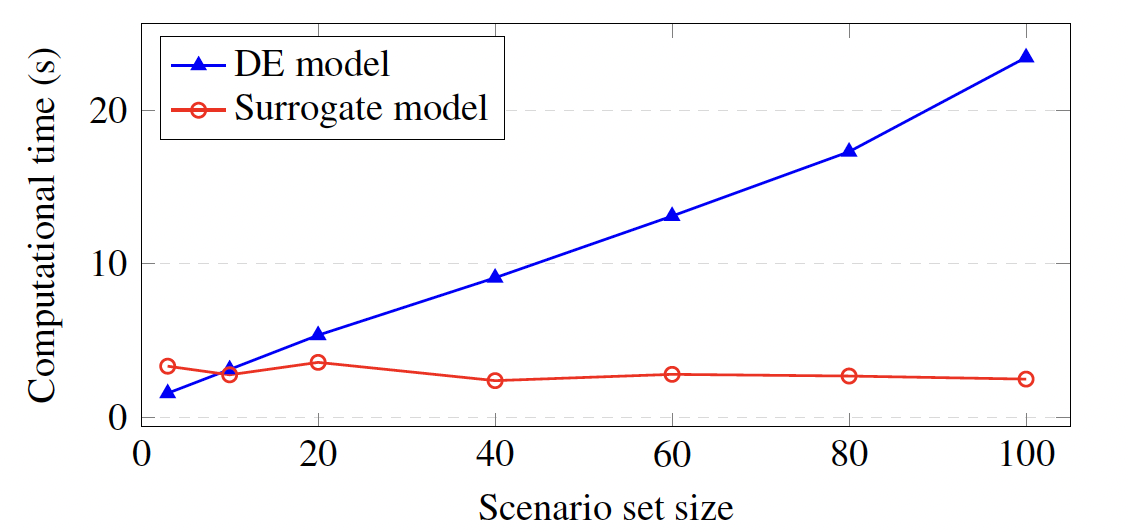}
        \caption{56 nodes - 50,000 observations}
        \label{fig:sub2}
    \end{subfigure}

    \vspace{1em}
    
    \begin{subfigure}[b]{0.46\textwidth}
        \centering
        \includegraphics[width=\textwidth]{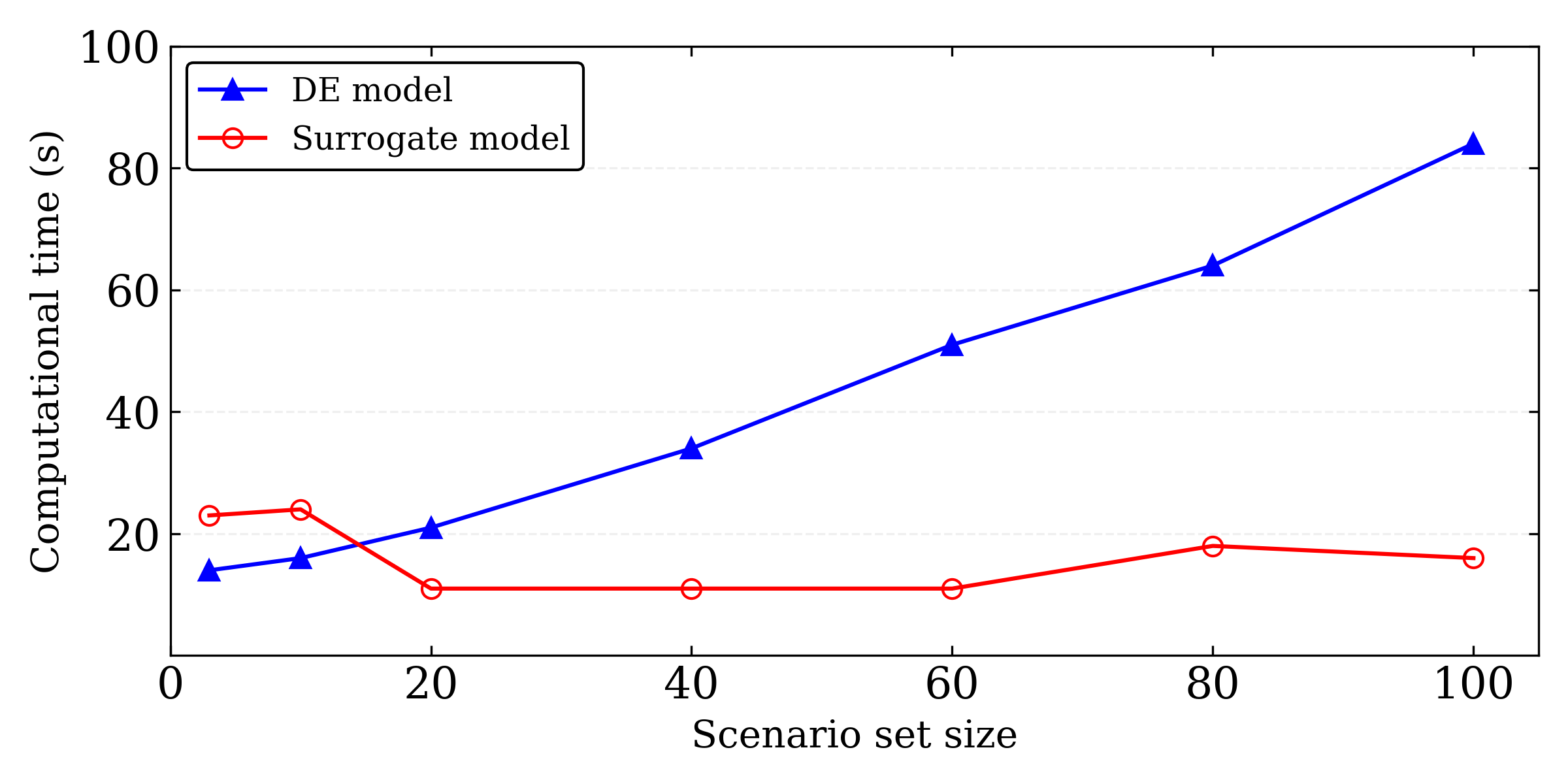}
        \caption{105 nodes - 100,000 observations}
        \label{fig:sub3}
    \end{subfigure}
    \hfill
    \begin{subfigure}[b]{0.49\textwidth}
        \centering
        \includegraphics[width=\textwidth]{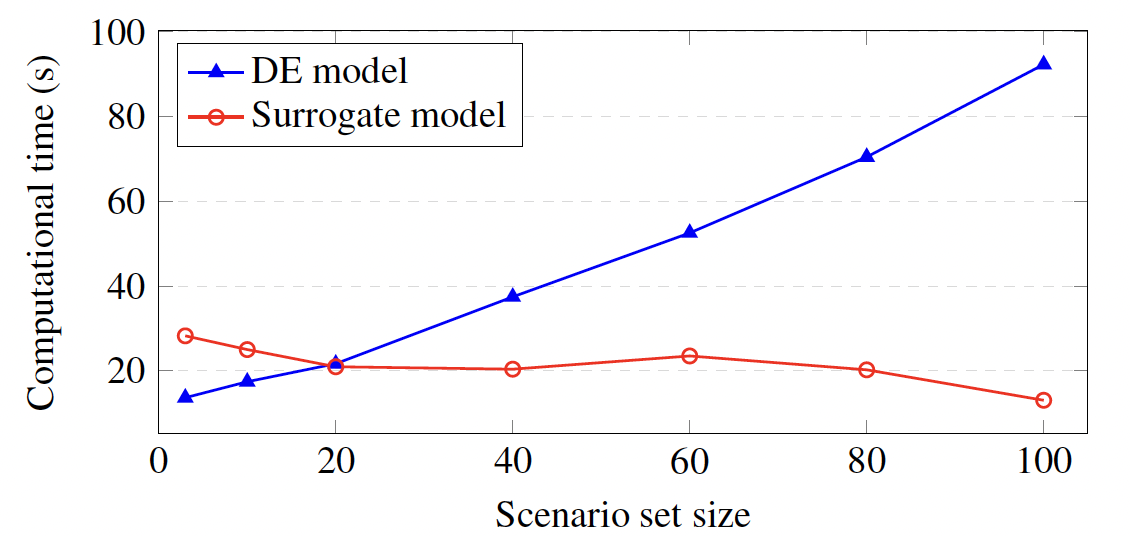}
        \caption{105 nodes - 200,000 observations}
        \label{fig:sub4}
    \end{subfigure}

    \caption{Computational time of NN surrogate and DE models versus scenario-set size}
    \label{fig:time}
\end{figure}

Figure~\ref{fig:time} reports the solving times for both the deterministic equivalent (DE model) and neural network surrogate (Surrogate model) approaches. The DE model shows an approximately linear increase in solving time with respect to the number of scenarios, as each additional scenario introduces a full replication of the recourse problem, thereby proportionally enlarging the MILP. It is worth noting that the flow-based DE model already represents a computationally efficient formulation compared to alternative DE approaches, such as subtour-elimination models \cite{Perboli17}, and therefore provides a meaningful benchmark.

In contrast, the Surrogate model maintains nearly constant solving times across all scenario set sizes. This behavior arises because the recourse problem is replaced by a compact surrogate model, whose size depends only on the fixed neural network architecture and is therefore independent of the number of scenarios.

This difference in scaling is highly beneficial for large-scale instances and extensive scenario sets, where the DE model incurs substantially higher computational costs. Consequently, the surrogate model achieves significant reductions in solving time while maintaining comparable solution quality and effectively preserving the stochastic characteristics of the original two-stage mpTSP.

\section{Conclusions}\label{sec:ccl}

In this work, we demonstrated that neural network models can be effectively integrated into a two-stage stochastic programming framework as surrogate models to generate high-quality solutions for the mpTSP under uncertainty, while significantly reducing computational time.
The experimental results indicate that surrogate models trained on relatively small datasets, even if showing slightly lower predictive accuracy, can produce solutions that are often more robust than those obtained from the Deterministic Equivalent model. As the computational effort of the surrogate-based optimization primarily depends on the neural network architecture, it is important to select a compact model that achieves an appropriate balance between predictive accuracy and computational efficiency.

Overall, the proposed neural network surrogate provides a scalable and generalizable approach by approximating the recourse function in a compact, differentiable, and solver-compatible form. Once trained, the surrogate can be efficiently reused for similar traffic conditions or incrementally updated with new data, thereby maintaining adaptability to evolving environments.

For future research, we plan to investigate the use of Graph Neural Networks to better exploit route topology and capture the combinatorial structure inherent in multi-path routing problems. In addition, we aim to develop surrogate models capable of generalizing across diverse problem instances, enabling a single pretrained model to provide accurate approximations without requiring full retraining. We also intend to extend the framework to multi-stage stochastic programming, where decisions and uncertainties evolve over multiple time periods. This extension will require the design of multi-step surrogate models and temporally coherent scenario generation mechanisms that reflect the dynamic structure of uncertainty. Such developments would allow the proposed framework to address large-scale, time-dependent network problems, supporting adaptive and high-quality decision-making in a wide range of real-world transportation applications.

\bmhead{Acknowledgments.} This work was supported by the Italian Ministry of University and Research through PRIN2020 project grant number 20207C8T9M.


\end{document}